\documentclass{amsart}
\usepackage{amstext}
\usepackage{amsthm}
\usepackage{amssymb}
\usepackage{esint}

\numberwithin{equation}{section}
\numberwithin{figure}{section}
\theoremstyle{plain}
\newtheorem{thm}{Theorem}
  \theoremstyle{definition}
  \newtheorem*{example*}{Example}
  \theoremstyle{remark}
  \newtheorem*{rem*}{Remark}

\allowdisplaybreaks
\usepackage[sort&compress,square,numbers]{natbib}
\usepackage{tikz}
\usetikzlibrary{matrix,positioning}

\makeatother

\begin{document}

\title{A note on Hermite polynomials}

\author{Taekyun Kim}
\address{Department of Mathematics, Kwangwoon University, Seoul 139-701, Republic
of Korea}
\email{tkkim@kw.ac.kr}

\author{Dae San Kim}
\address{Department of Mathematics, Sogang University, Seoul 121-742, Republic
of Korea}
\email{dskim@sogang.ac.kr}

\begin{abstract}
In this paper, we consider linear differential equations satisfied by the generating function for Hermite polynomials and derive some new identities involving those polynomials.
\end{abstract}



\maketitle
\global\long\def\relphantom#1{\mathrel{\phantom{{#1}}}}

\section{Introduction}

The Hermite polynomials form a Sheffer sequence and are given by
the generating function  
\begin{equation}
e^{2xt-t^{2}}=\sum_{n=0}^{\infty}\frac{H_{n}\left(x\right)}{n!}t^{n},\quad\left(\text{see \cite{key-1,key-2,key-3,key-4,key-5,key-6,key-7,key-8,key-9,key-10,key-11}}\right).\label{eq:1}
\end{equation}

By using Taylor series, we get 
\begin{align*}
H_{n}\left(x\right) & =\left[\left(\frac{\partial}{\partial t}\right)^{n}e^{\left(2xt-t^{2}\right)}\right]_{t=0}\\
 & =\left[e^{x^{2}}\left(\frac{\partial}{\partial t}\right)^{n}e^{-\left(x-t\right)^{2}}\right]_{t=0}\\
 & =\left(-1\right)^{n}e^{x^{2}}\left[\left(\frac{\partial}{\partial x}\right)^{n}e^{-\left(x-t\right)^{2}}\right]_{t=0}\\
 & =\left(-1\right)^{n}e^{x^{2}}\frac{d^{n}}{dx^{n}}e^{-x^{2}},\quad\left(n\ge0\right),\quad\left(\text{see \cite{key-1,key-2,key-3,key-4,key-5,key-6,key-7,key-8,key-9,key-10,key-11,key-12,key-13,key-14,key-15,key-16}}\right).
\end{align*}

The Hermite polynomials can be represented by the Contour integral
as follows: 
\begin{equation}
H_{n}\left(z\right)=\frac{n!}{2\pi i}\oint e^{-t^{2}+2tz}t^{-n-1}dt,\label{eq:2}
\end{equation}
where the Contour encloses the origin and is traversed in a counterclockwise
direction (see \cite{key-2,key-8,key-11,key-14}).

The probabilists' Hermite polynomials are given by the generating function
\begin{align}
H_{n}^{*}\left(x\right) & =\left(-1\right)^{n}e^{\frac{x^{2}}{2}}\frac{d^{n}}{dx^{n}}e^{-\frac{x^{2}}{2}}\label{eq:3}\\
 & =\left(x-\frac{d}{dx}\right)^{n}\cdot1,\quad\left(\text{see \cite{key-9}}\right).\nonumber 
\end{align}

The physicists' Hermite polynomials are also given by 
\begin{align}
H_{n}\left(x\right) & =\left(-1\right)^{n}e^{x^{2}}\frac{d^{n}}{dx^{n}}e^{-x^{2}}\label{eq:4}\\
 & =\left(2x-\frac{d}{dx}\right)^{n}\cdot1\quad\left(\text{see \cite{key-20}}\right).\nonumber 
\end{align}

Thus, by (\ref{eq:3}) and (\ref{eq:4}), we get 
\begin{equation}
H_{n}\left(x\right)=2^{\frac{n}{2}}H_{n}^{*}\left(\sqrt{2}x\right),\quad H_{n}^{*}\left(x\right)=2^{-\frac{n}{2}}H_{n}\left(\frac{x}{\sqrt{2}}\right),\label{eq:5}
\end{equation}
where $n\ge0$ (see \cite{key-12,key-13,key-14,key-15,key-16}). 

The first several Hermite polynomials are $H_{0}\left(x\right)=1$,
$H_{1}\left(x\right)=2x$, $H_{2}\left(x\right)=4x^{2}-2$, $H_{3}\left(x\right)=8x^{2}-12x$,
$H_{4}\left(x\right)=16x^{4}-48x^{2}+12$, $H_{5}\left(x\right)=32x^{5}-160x^{3}+120x$,
$H_{6}\left(x\right)=64x^{6}-480x^{4}+720x^{2}-120$, ...

The probabilists' Hermite polynomials are solutions of the differential
equation: 
\[
\left(e^{-\frac{x^{2}}{2}}u^{\prime}\right)^{\prime}+\lambda e^{-\frac{1}{2}x^{2}}u=0,
\]
where $\lambda$ is a constant, with the boundary conditions that $u$
should be polynomially bounded at infinity. 

The generating function of the probabilists' Hermite polynomials is
given by 
\begin{equation}
e^{xt-\frac{t^{2}}{2}}=\sum_{n=0}^{\infty}H_{n}^{*}\left(x\right)\frac{t^{n}}{n!},\quad\left(\text{see \cite{key-13,key-15,key-16}}\right).\label{eq:6}
\end{equation}

The Hermite polynomials $H_{n}^{\left(\nu\right)}\left(x\right)$
of variance $\nu$ form an Appell sequence and are defined by the generating
function  
\begin{equation}
\sum_{k=0}^{\infty}\frac{H_{k}^{\left(\nu\right)}\left(x\right)}{k!}t^{k}=e^{xt-\frac{\nu t^{2}}{2}},\quad\left(\text{see \cite{key-13}}\right).\label{eq:7}
\end{equation}

Thus, by (\ref{eq:7}), we get 
\begin{equation}
x^{2m+1}=\sum_{l=0}^{m}\binom{2m+1}{2l+1}\frac{\left(2m-2l\right)!}{\left(m-l\right)!}\left(\frac{\nu}{2}\right)^{m-l}H_{2l+1}^{\left(\nu\right)}\left(x\right),\label{eq:8}
\end{equation}
and 
\begin{equation}
x^{2m}=\sum_{l=0}^{m}\binom{2m}{2l}\frac{\left(2m-2l\right)!}{\left(m-l\right)!}\left(\frac{\nu}{2}\right)^{m-l}H_{2l}^{\left(\nu\right)}\left(x\right),\quad\left(\text{see \cite{key-13}}\right).\label{eq:9}
\end{equation}

The Hermite polynomials have been studied in probability, combinatorics,
numerical analysis, finite element methods, physics and system theory
(see \cite{key-1,key-2,key-3,key-4,key-5,key-6,key-7,key-8,key-9,key-10,key-11,key-12,key-13,key-14,key-15,key-16}). 

Recently, Kim has studied nonlinear differential equations arising
from Frobenius-Euler numbers and polynomials.

In this paper, we consider linear differential equations arising from
Hermite polynomials of variance $\nu$ and give some new and explicit
identities for those polynomials.

\section{Hermite polynomials of variance $\nu$}

Let 
\begin{equation}
F=F\left(t:x,\nu\right)=e^{xt-\frac{\nu t^{2}}{2}}.\label{eq:10}
\end{equation}

From (\ref{eq:10}), we note that 
\begin{align}
F^{\left(1\right)} & =\frac{d}{dt}F\left(t:x,\nu\right)\label{eq:11}\\
 & =\left(x-\nu t\right)e^{xt-\frac{\nu t^{2}}{2}}\nonumber \\
 & =\left(x-\nu t\right)F,\nonumber 
\end{align}
\begin{align}
F^{\left(2\right)} & =\frac{d}{dt}F^{\left(1\right)}=\left(-\nu+\left(x-t\nu\right)^{2}\right)F,\label{eq:12}\\
F^{\left(3\right)} & =\frac{d}{dt}F^{\left(2\right)}=\left(-3\nu\left(x-\nu t\right)+\left(x-\nu t\right)^{3}\right)F,\label{eq:13}
\end{align}
and 
\begin{equation}
F^{\left(4\right)}=\frac{d}{dt}F^{\left(3\right)}=\left(3\nu^{2}-6\nu\left(x-\nu t\right)^{2}+\left(x-\nu t\right)^{4}\right)F.\label{eq:14}
\end{equation}

Continuing this process, we set 
\begin{align}
F^{\left(N\right)} & =\left(\frac{d}{dt}\right)^{N}F\left(t:x,\nu\right)\label{eq:15}\\
 & =\left(\sum_{i=0}^{N}a_{i}\left(N,\nu\right)\left(x-\nu t\right)^{i}\right)F,\nonumber 
\end{align}
where $N\in\mathbb{N}\cup\left\{ 0\right\} $.

From (\ref{eq:15}), we have 
\begin{align}
F^{\left(N+1\right)} & =\frac{d}{dt}F^{\left(N\right)}\label{eq:16}\\
 & =\sum_{i=0}^{N}a_{i}\left(N,\nu\right)i\left(x-\nu t\right)^{i-1}\left(-\nu\right)F\nonumber \\
 & \relphantom =+\sum_{i=0}^{N}a_{i}\left(N,\nu\right)\left(x-\nu t\right)^{i}F^{\left(1\right)}.\nonumber 
\end{align}

By (\ref{eq:11}) and (\ref{eq:16}), we easily get 
\begin{align}
F^{\left(N+1\right)} & =\left\{ -\nu a_{1}\left(N,\nu\right)+a_{N}\left(N,\nu\right)\left(x-\nu t\right)^{N+1}+a_{N-1}\left(N,\nu\right)\left(x-\nu t\right)^{N}\right.\label{eq:17}\\
 & \relphantom =+\left.\sum_{i=1}^{N-1}\left(-\left(i+1\right)\nu a_{i+1}\left(N,\nu\right)+a_{i-1}\left(N,\nu\right)\right)\left(x-\nu t\right)^{i}\right\} F.\nonumber 
\end{align}

By replacing $N$ by $\left(N+1\right)$ in (\ref{eq:15}), we get
\begin{equation}
F^{\left(N+1\right)}=\left(\sum_{i=0}^{N+1}a_{i}\left(N+1,\nu\right)\left(x-\nu t\right)^{i}\right)F.\label{eq:18}
\end{equation}

From (\ref{eq:17}) and (\ref{eq:18}), we can derive the following
equations: 
\begin{align}
a_{0}\left(N+1,\nu\right) & =-\nu a_{1}\left(N,\nu\right),\label{eq:19}\\
a_{N}\left(N+1,\nu\right) & =a_{N-1}\left(N,\nu\right),\label{eq:20}\\
a_{N+1}\left(N+1,\nu\right) & =a_{N}\left(N,\nu\right)\label{eq:21}
\end{align}
and 
\begin{equation}
a_{i}\left(N+1,\nu\right)=-\left(i+1\right)\nu a_{i+1}\left(N,\nu\right)+a_{i-1}\left(N,\nu\right),\label{eq:22}
\end{equation}
where $1\le i\le N-1$.

It is not difficult to show that 
\begin{equation}
F=F^{\left(0\right)}=a_{0}\left(0,\nu\right)F.\label{eq:23}
\end{equation}

Thus, by (\ref{eq:23}), we get 
\begin{equation}
a_{0}\left(0,\nu\right)=1.\label{eq:24}
\end{equation}

From (\ref{eq:11}) and (\ref{eq:15}), we note that 
\begin{equation}
\left(x-\nu t\right)F=F^{\left(1\right)}=\left(a_{0}\left(1,\nu\right)+a_{1}\left(1,\nu\right)\left(x-\nu t\right)\right)F.\label{eq:25}
\end{equation}

Thus, by comparing the coefficients on both sides of (\ref{eq:25}),
we get 
\begin{equation}
a_{0}\left(1,\nu\right)=0,\quad a_{1}\left(1,\nu\right)=1.\label{eq:26}
\end{equation}

From (\ref{eq:20}), (\ref{eq:21}), (\ref{eq:24}) and (\ref{eq:26}),
we have 
\begin{equation}
a_{N}\left(N+1,\nu\right)=a_{N-1}\left(N,\nu\right)=\cdots=a_{0}\left(1,\nu\right)=0,\label{eq:27}
\end{equation}
and 
\begin{equation}
a_{N+1}\left(N+1,\nu\right)=a_{N}\left(N,\nu\right)=\cdots=a_{1}\left(1,\nu\right)=1.\label{eq:28}
\end{equation}

Therefore, we obtain the following theorem.
\begin{thm}
\label{thm:1} The linear differential equations
\begin{align*}
F^{\left(N\right)} & =\left(\frac{d}{dt}\right)^{N}F\left(t:x,\nu\right)\\
 & =\left(\sum_{i=0}^{N}a_{i}\left(N,\nu\right)\left(x-\nu t\right)^{i}\right)F,\quad\left(N\in\mathbb{N}\cup\left\{ 0\right\} \right)
\end{align*}
has a solution $F=F\left(t:x,\nu\right)=e^{xt-\frac{\nu t^{2}}{2}}$,
where 
\begin{align*}
a_{0}\left(N,\nu\right) & =-\nu a_{1}\left(N-1,\nu\right),\\
a_{N-1}\left(N,\nu\right) & =a_{N-2}\left(N-1,\nu\right)=\cdots=a_{1}\left(2,\nu\right)=a_{0}\left(1,\nu\right)=0,\\
a_{N}\left(N,\nu\right) & =a_{N-1}\left(N-1,\nu\right)=\cdots=a_{1}\left(1,\nu\right)=a_{0}\left(0,\nu\right)=1,
\end{align*}
and 
\[
a_{i}\left(N,\nu\right)=-\left(i+1\right)\nu a_{i+1}\left(N-1,\nu\right)+a_{i-1}\left(N-1,\nu\right),\quad\left(1\le i\le N-2\right).
\]
\end{thm}
\begin{example*}
$\:$
\begin{enumerate}
\item $N=3$, $i=1$. By (\ref{eq:22}), we get 
\begin{align*}
a_{1}\left(3,\nu\right) & =-2\nu a_{2}\left(2,\nu\right)+a_{0}\left(2,\nu\right)\\
 & =-2\nu-\nu=-3\nu.
\end{align*}

\item $N=4$, $1\le i\le2$. By (\ref{eq:22}), we have 
\[
a_{1}\left(4,\nu\right)=0,\quad a_{2}\left(4,\nu\right)=-6\nu.
\]

\item $N=5$, $1\le i\le3$. By (\ref{eq:22}), we get 
\[
a_{1}\left(5,\nu\right)=15\nu^{2},\quad a_{2}\left(5,\nu\right)=0,\quad a_{3}\left(5,\nu\right)=-10\nu.
\]

\item $N=6$, $1\le i\le4$. From (\ref{eq:22}), we have 
\[
a_{1}\left(6,\nu\right)=0,\quad a_{2}\left(6,\nu\right)=45\nu^{2},\quad a_{3}\left(6,\nu\right)=0,\quad a_{4}\left(6,\nu\right)=-15\nu.
\]

\end{enumerate}
\end{example*}
Thus, we obtain the following result.

\begin{rem*}
$\,$
The matrix $\left(a_{i}\left(j,\nu\right)\right)_{0\le i,j\le6}$
is given by 

\begin{equation*}
 \begin{tikzpicture}[baseline=(current  bounding  box.west)]
  \matrix (mymatrix) [matrix of math nodes,left delimiter={[},right
delimiter={]}]
  {
   1 & 0 & -\nu & 0 & 3\nu^{2} & 0 & -15\nu^{3}\\
 & 1 & 0 & -3\nu & 0 & 15\nu^{2} & 0\\
 &  & 1 & 0 & -6\nu & 0 & 45\nu^{2}\\
 &  &  & 1 & 0 & -10\nu & 0\\
 &  &  &  & 1 & 0 & -15\nu\\
 &  &  &  &  & 1 & 0\\
 &  &  &  &  &  & 1\\
  };
\node[xshift=-74pt,yshift=59pt] {$0$};
\node[xshift=-62pt,yshift=59pt] {$1$};
\node[xshift=-45pt,yshift=59pt] {$2$};
\node[xshift=-24pt,yshift=59pt] {$3$};
\node[xshift=1pt,yshift=59pt] {$4$};
\node[xshift=28pt,yshift=59pt] {$5$};
\node[xshift=62pt,yshift=59pt] {$6$};
\node[xshift=-95pt,yshift=40pt] {$0$};
\node[xshift=-95pt,yshift=26pt] {$1$};
\node[xshift=-95pt,yshift=11pt] {$2$};
\node[xshift=-95pt,yshift=-4pt] {$3$};
\node[xshift=-95pt,yshift=-18pt] {$4$};
\node[xshift=-95pt,yshift=-32pt] {$5$};
\node[xshift=-95pt,yshift=-46pt] {$6$};
\node[xshift=-40pt,yshift=-30pt] {\LARGE$0$};
\end{tikzpicture}.
\end{equation*}

\end{rem*}

From (\ref{eq:7}), we note that 
\begin{align}
F & =F\left(t:x,\nu\right)=e^{xt-\frac{\nu t^{2}}{2}}\label{eq:29}\\
 & =\sum_{k=0}^{\infty}H_{k}^{\left(\nu\right)}\left(x\right)\frac{t^{k}}{k!}.\nonumber 
\end{align}

Thus, by (\ref{eq:29}), we get 
\begin{align}
F^{\left(N\right)} & =\left(\frac{d}{dt}\right)^{N}F\left(t:x,\nu\right)\label{eq:30}\\
 & =\sum_{k=N}^{\infty}H_{k}^{\left(\nu\right)}\left(x\right)\left(k\right)_{N}\frac{t^{k-N}}{k!}\nonumber \\
 & =\sum_{k=0}^{\infty}H_{k+N}^{\left(\nu\right)}\left(x\right)\left(k+N\right)_{N}\frac{t^{k}}{\left(n+k\right)!}\nonumber \\
 & =\sum_{k=0}^{\infty}H_{k+N}^{\left(\nu\right)}\left(x\right)\frac{t^{k}}{k!}.\nonumber 
\end{align}

By Theorem \ref{thm:1}, we easily get 
\begin{align}
 & F^{\left(N\right)}\label{eq:31}\\
 & =\left(\sum_{i=0}^{N}a_{i}\left(N,\nu\right)\left(x-\nu t\right)^{i}\right)F\nonumber \\
 & =\sum_{i=0}^{N}a_{i}\left(N,\nu\right)\sum_{m=0}^{\infty}\left(i\right)_{m}x^{i-m}\left(-\nu\right)^{m}\frac{t^{m}}{m!}\sum_{l=0}^{\infty}H_{l}^{\left(\nu\right)}\left(x\right)\frac{t^{l}}{l!}\nonumber \\
 & =\sum_{k=0}^{\infty}\left\{ \sum_{i=0}^{N}a_{i}\left(N,\nu\right)\sum_{l=0}^{k}\binom{k}{l}\left(i\right)_{k-l}\left(-\nu\right)^{k-l}x^{i+l-k}H_{l}^{\left(\nu\right)}\left(x\right)\right\} \frac{t^{k}}{k!}\nonumber \\
 & =\sum_{k=0}^{\infty}\left\{ \sum_{i=0}^{N}a_{i}\left(N,\nu\right)\sum_{l=\max\left\{ 0,k-i\right\} }^{k}\binom{k}{l}\left(i\right)_{k-l}\left(-\nu\right)^{k-l}x^{i+l-k}H_{l}^{\left(\nu\right)}\left(x\right)\right\} \frac{t^{k}}{k!}.\nonumber 
\end{align}

Therefore, by (\ref{eq:30}) and (\ref{eq:31}), we obtain the following
theorem.
\begin{thm}
\label{thm:2} For $k,N\in\mathbb{N}\cup\left\{ 0\right\} $, we have
\begin{align*}
 & H_{k+N}^{\left(\nu\right)}\left(x\right)\\
 & =\sum_{i=0}^{N}a_{i}\left(N,\nu\right)\sum_{l=\max\left\{ 0,k-i\right\} }^{k}\binom{k}{l}\left(i\right)_{k-l}\left(-\nu\right)^{k-l}x^{i+l-k}H_{l}^{\left(\nu\right)}\left(x\right).
\end{align*}

\end{thm}
It is easy to show that 
\begin{equation}
H_{k+1}^{\left(\nu\right)}\left(x\right)=\left(x-\nu\frac{\partial}{\partial x}\right)H_{k}^{\left(\nu\right)}\left(x\right).\label{eq:32}
\end{equation}

Thus, by (\ref{eq:32}), we have 

\begin{equation}
H_{k+N}^{\left(\nu\right)}\left(x\right)=\left(x-\nu\frac{\partial}{\partial x}\right)^{N}H_{k}^{\left(\nu\right)}\left(x\right),\quad\left(N\in\mathbb{N}\cup\left\{ 0\right\} \right).\label{eq:33}
\end{equation}

From Theorem \ref{thm:2}, we note that 
\begin{align}\label{eq:33-1}
&\relphantom{=}\left(x-\nu\frac{\partial}{\partial x}\right)^{N}H_{k}^{\left(\nu\right)}\left(x\right)\\
&=\sum_{i=0}^{N}a_{i}\left(N,\nu\right)\sum_{l=\max\left\{ 0,k-i\right\} }^{k}\binom{k}{l}\left(i\right)_{k-l}\left(-\nu\right)^{k-l}x^{i+l-k}H_{l}^{\left(\nu\right)}\left(x\right),\nonumber
\end{align}
 where $\frac{\partial}{\partial x}x-x\frac{\partial}{\partial x}=\text{identity}$. 

Now, we observe explicit determination of $a_{i}\left(j,\nu\right)$. 

From (\ref{eq:21}) and (\ref{eq:22}), we can derive the following
equations: 
\begin{align}
a_{N}\left(N,\nu\right) & =1,\label{eq:34}
\end{align}
\begin{align}
a_{N-2}\left(N,\nu\right) & =-\left(N-1\right)\nu a_{N-1}\left(N-1,\nu\right)+a_{N-3}\left(N-1,\nu\right)\label{eq:35}\\
 & =-\left(N-1\right)\nu a_{N-1}\left(N-1,\nu\right)-\left(N-2\right)\nu a_{N-2}\left(N-2,\nu\right)\nonumber\\
 &\relphantom{=}+a_{N-4}\left(N-2,\nu\right)\nonumber \\
 & \vdots\nonumber \\
 & =-\left(N-1\right)\nu a_{N-1}\left(N-1,\nu\right)-\left(N-2\right)\nu a_{N-2}\left(N-2,\nu\right)\nonumber\\
 &\relphantom{=}-\cdots-2\nu a_{2}\left(2,\nu\right)+a_{0}\left(2,\nu\right)\nonumber \\
 & =-\left(N-1\right)\nu a_{N-1}\left(N-1,\nu\right)-\left(N-2\right)\nu a_{N-2}\left(N-2,\nu\right)\nonumber\\
 &\relphantom{=}-\cdots-2\nu a_{2}\left(2,\nu\right)-\nu a_{1}\left(1,\nu\right)\nonumber \\
 & =-\nu\sum_{i=1}^{N-1}ia_{i}\left(i,\nu\right),\nonumber 
\end{align}

\begin{align}
a_{N-4}\left(N,\nu\right) & =-\left(N-3\right)\nu a_{N-3}\left(N-1,\nu\right)+a_{N-5}\left(N-1,\nu\right)\label{eq:36}\\
 & =-\left(N-3\right)\nu a_{N-3}\left(N-1,\nu\right)-\left(N-4\right)\nu a_{N-4}\left(N-2,\nu\right)\nonumber\\
 &\relphantom{=}+a_{N-6}\left(N-2,\nu\right)\nonumber \\
 & \vdots\nonumber \\
 & =-\left(N-3\right)\nu a_{N-3}\left(N-1,\nu\right)-\left(N-4\right)\nu a_{N-4}\left(N-2,\nu\right)\nonumber\\
 &\relphantom{=}-\cdots-2\nu a_{2}\left(4,\nu\right)+a_{0}\left(4,\nu\right)\nonumber \\
 & =-\left(N-3\right)\nu a_{N-3}\left(N-1,\nu\right)-\left(N-4\right)\nu a_{N-4}\left(N-2,\nu\right)\nonumber\\
 &\relphantom{=}-\cdots-2\nu a_{2}\left(4,\nu\right)-\nu a_{1}\left(3,\nu\right)\nonumber \\
 & =-\nu\sum_{i=0}^{N-3}ia_{i}\left(i+2,\nu\right),\nonumber 
\end{align}
and 
\begin{align}
a_{N-6}\left(N,\nu\right) & =-\left(N-5\right)\nu a_{N-5}\left(N-1,\nu\right)+a_{N-7}\left(N-1,\nu\right)\label{eq:37}\\
 & =-\left(N-5\right)\nu a_{N-5}\left(N-1,\nu\right)-\left(N-6\right)\nu a_{N-6}\left(N-2,\nu\right)\nonumber\\
 &\relphantom{=}+a_{N-8}\left(N-2,\nu\right)\nonumber \\
 & \vdots\nonumber \\
 & =-\left(N-5\right)\nu a_{N-5}\left(N-1,\nu\right)-\left(N-6\right)\nu a_{N-6}\left(N-2,\nu\right)\nonumber\\
 &\relphantom{=}-\cdots-2\nu a_{2}\left(6,\nu\right)-\nu a_{1}\left(5,\nu\right)\nonumber \\
 & =-\nu\sum_{i=1}^{N-5}ia_{i}\left(i+4,\nu\right).\nonumber 
\end{align}

Continuing in this fashion, for $l$ with  $1\le l\le\left[\frac{N-1}{2}\right]$,
\begin{equation}
a_{N-2l}\left(N,\nu\right)=-\nu\sum_{i=1}^{N-2l+1}ia_{i}\left(i+2l-2,\nu\right).\label{eq:38}
\end{equation}

By (\ref{eq:34}), (\ref{eq:35}), (\ref{eq:36}), (\ref{eq:37})
and (\ref{eq:38}), we get 
\begin{align}
a_{N-2}\left(N,\nu\right) & =-\nu\sum_{i_{1}=1}^{N-1}i_{1},\label{eq:39}\\
a_{N-4}\left(N,\nu\right) & =-\nu\sum_{i_{2}=1}^{N-3}i_{2}a_{i_{2}}\left(i_{2}+2,\nu\right)\label{eq:40}\\
 & =\left(-\nu\right)^{2}\sum_{i_{2}=1}^{N-3}\sum_{i_{1}=1}^{i_{2}+1}i_{2}i_{1},\nonumber \\
a_{N-6}\left(N,\nu\right) & =-\nu\sum_{i_{3}=1}^{N-5}i_{3}a_{i_{3}}\left(i_{3}+4,\nu\right)\label{eq:41}\\
 & =\left(-\nu\right)^{3}\sum_{i_{3}=1}^{N-5}\sum_{i_{2}=1}^{i_{3}+1}\sum_{i_{1}=1}^{i_{2}+1}i_{3}i_{2}i_{1},\nonumber 
\end{align}
and 
\begin{equation}
a_{N-2l}\left(N,\nu\right)=\left(-\nu\right)^{l}\sum_{i_{l}=1}^{N-2l+1}\sum_{i_{l-1}=1}^{i_{l}+1}\cdots\sum_{i_{1}=1}^{i_{2}+1}i_{l}\cdot i_{l-1}\cdots i_{1},\label{eq:42}
\end{equation}
where $1\le l\le\left[\frac{N-1}{2}\right]$.

By (\ref{eq:20}) and (\ref{eq:22}), we easily get 
\begin{align}
a_{N-1}\left(N,\nu\right) & =a_{N-2}\left(N-1,\nu\right)=a_{N-3}\left(N-2,\nu\right)=\cdots=a_{0}\left(1,\nu\right)=0,\label{eq:43}\\
a_{N-3}\left(N,\nu\right) & =-\left(N-2\right)\nu a_{N-2}\left(N-1,\nu\right)+a_{N-4}\left(N-1,\nu\right)\label{eq:44}\\
 & =a_{N-4}\left(N-1,\nu\right)\nonumber \\
 & \vdots\nonumber \\
 & =a_{0}\left(3,\nu\right)=-\nu a_{1}\left(2,\nu\right)=-\nu a_{0}\left(1,\nu\right)=0,\nonumber \\
a_{N-5}\left(N,\nu\right) & =-\left(N-4\right)\nu a_{N-4}\left(N-1,\nu\right)+a_{N-6}\left(N-1,\nu\right)=a_{N-6}\left(N-1,\nu\right)\label{eq:45}\\
 & \vdots\nonumber \\
 & =a_{0}\left(5,\nu\right)=-\nu a_{1}\left(4,\nu\right)=0,\nonumber \\
a_{N-7}\left(N,\nu\right) & =-\left(N-6\right)\nu a_{N-6}\left(N-1,\nu\right)+a_{N-8}\left(N-1,\nu\right)\label{eq:46}\\
 & \vdots\nonumber \\
 & =a_{0}\left(7,\nu\right)=-\nu a_{1}\left(6,\nu\right)=0,\nonumber 
\end{align}
and 
\begin{equation}
a_{N-\left(2l-1\right)}\left(N,\nu\right)=0,\quad\left(1\le l\le\left[\frac{N}{2}\right]\right).\label{eq:47}
\end{equation}

Therefore, we obtain the following theorem.
\begin{thm}
\label{thm:3} For $N\in\mathbb{N}\cup\left\{ 0\right\} $, we have
\[
a_{N-2l}\left(N,\nu\right)=\left(-\nu\right)^{l}\sum_{i_{l}=1}^{N-2l+1}\sum_{i_{l-1}=1}^{i_{l}+1}\cdots\sum_{i_{1}=1}^{i_{2}+1}i_{l}i_{l-1}\cdots i_{1},
\]
where $1\le l\le\left[\frac{N-1}{2}\right]$.

Also, 
\[
a_{N-\left(2l-1\right)}\left(N,\nu\right)=0,\quad\text{if }1\le l\le\left[\frac{N}{2}\right].
\]

\end{thm}
\bibliographystyle{amsplain}
\nocite{*}
\providecommand{\bysame}{\leavevmode\hbox to3em{\hrulefill}\thinspace}
\providecommand{\MR}{\relax\ifhmode\unskip\space\fi MR }
\providecommand{\MRhref}[2]{%
  \href{http://www.ams.org/mathscinet-getitem?mr=#1}{#2}
}
\providecommand{\href}[2]{#2}

\end{document}